\documentclass[12pt]{amsart}
\usepackage{amssymb,amsmath,amsfonts,latexsym, amscd}
\usepackage{bm}
\newcommand{\newsection}[1]{\setcounter{equation}{0} \section{#1}}
\def\V{{\underline{V}}}

\newcommand{\bea}{\begin{eqnarray}}
\newcommand{\eea}{\end{eqnarray}}

\newcommand{\cld}{\mathcal{D}}
\newcommand{\cle}{\mathcal{E}}

\newcommand{\clh}{\mathcal{H}}
\newcommand{\clk}{\mathcal{K}}
\newcommand{\cll}{\mathcal{L}}
\newcommand{\clm}{\mathcal{M}}
\newcommand{\cln}{\mathcal{N}}

\newcommand{\raro}{\rightarrow}

\def \qed {\hfill \vrule height6pt width 6pt depth 0pt}
\def\textmatrix#1&#2\\#3&#4\\{\bigl({#1 \atop #3}\ {#2 \atop #4}\bigr)}
\def\dispmatrix#1&#2\\#3&#4\\{\left({#1 \atop #3}\ {#2 \atop #4}\right)}
\newcommand{\be}{\begin{equation}}
\newcommand{\ee}{\end{equation}}
\newcommand{\ben}{\begin{eqnarray*}}
\newcommand{\een}{\end{eqnarray*}}

\newcommand{\NI}{\noindent}

\newcommand{\bi}{\begin{itemize}}
\newcommand{\ei}{\end{itemize}}

\def\5{{5\superprime}}

\newtheorem{Theorem}{\sc Theorem}[section]
\newtheorem{Lemma}[Theorem]{\sc Lemma}
\newtheorem{Proposition}[Theorem]{\sc Proposition}
\newtheorem{Corollary}[Theorem]{\sc Corollary}
\newtheorem{Definition}[Theorem]{\sc Definition}
\newtheorem{Example}[Theorem]{\sc Example}

\newtheorem{Note}[Theorem]{\sc Note}
\newtheorem{Question}{\sc Question}
\newtheorem{ass}[Theorem]{\sc Assumption}

\theoremstyle{remark}
\newtheorem{Remark}[Theorem]{\sc Remark}

\newcommand{\bt}{\begin{Theorem}}
\def\beginlem{\begin{Lemma}}
\def\beginprop{\begin{Proposition}}
\def\begincor{\begin{Corollary}}
\def\begindef{\begin{Definition}}
\def\beginexamp{\begin{Example}}
\def\beginrem{\begin{Remark}}
\def\beginq{\begin{Question}}
\def\beginass{\begin{ass}}
\def\beginnote{\begin{Note}}
\newcommand{\et}{\end{Theorem}}
\def\endlem{\end{Lemma}}
\def\endprop{\end{Proposition}}
\def\endcor{\end{Corollary}}
\def\enddef{\end{Definition}}
\def\endexamp{\end{Example}}
\def\endrem{\end{Remark}}
\def\endq{\end{Question}}
\def\endass{\end{ass}}
\def\endnote{\end{Note}}

\begin{document}

\title[Contractions with Polynomial characteristic functions]{Contractions with Polynomial characteristic functions II. Analytic approach}

\dedicatory{This paper is dedicated to our lifelong friend Ron
Douglas on the occasion of his upcoming 80th birthday}

\author[Foias]{Ciprian Foias}
\address[Ciprian Foias]{Department of Mathematics, Texas A\&M University, College Station, Texas 77843, USA}

\author[Pearcy]{Carl Pearcy}
\address[Carl Pearcy]{Department of Mathematics, Texas A\&M University, College Station, Texas 77843, USA}
\email{cpearcy@math.tamu.edu}

\author[Sarkar]{Jaydeb Sarkar}
\address[Jaydeb Sarkar]{Indian Statistical Institute, Statistics and Mathematics Unit, 8th Mile, Mysore Road, Bangalore, 560059, India}
\email{jay@isibang.ac.in, jaydeb@gmail.com}



\begin{abstract}
The simplest and most natural examples of completely nonunitary
contractions on separable complex Hilbert spaces which have
polynomial characteristic functions are the nilpotent operators. The
main purpose of this paper is to prove the following theorem: Let
$T$ be a completely nonunitary contraction on a Hilbert space
$\mathcal{H}$. If the characteristic function $\Theta_T$ of $T$ is a
polynomial of degree $m$, then there exist a Hilbert space
$\mathcal{M}$, a nilpotent operator $N$ of order $m$, a coisometry
$V_1 \in \mathcal{L}(\overline{\mbox{ran}} (I - N N^*) \oplus
\mathcal{M}, \overline{\mbox{ran}} (I - T T^*))$, and an isometry
$V_2 \in \mathcal{L}(\overline{\mbox{ran}} (I - T^* T),
\overline{\mbox{ran}} (I - N^* N) \oplus \mathcal{M})$, such that
\[
\Theta_T = V_1 \begin{bmatrix} \Theta_N & 0
\\ 0 & I_{\mathcal{M}} \end{bmatrix} V_2.
\]
\end{abstract}


\maketitle
\newsection{Introduction}

This is a sequel to our paper \cite{FS}, where we identified the
structure of the completely nonunitary contractions on a Hilbert
space that have a polynomial characteristic function. Namely, we
proved that the characteristic function $\Theta_T$ of a completely
nonunitary contraction $T$ on a separable, infinite dimensional,
complex Hilbert space $\clh$ is a polynomial if and only if there
exist three closed subspaces $\clh_1, \clh_0, \clh_{-1}$ of $\clh$
with $\clh = \clh_1 \oplus \clh_0 \oplus \clh_{-1}$, a pure isometry
$S$ in $\cll(\clh_1)$, a nilpotent $N$ in $\cll(\clh_0)$, and a pure
co-isometry $C$ in $\cll(\clh_{-1})$, such that $T$ has the matrix
representation
\[
T = \begin{bmatrix} S & * & *\\0 & N & *\\0& 0& C
\end{bmatrix}.
\]

Moreover, the multiplicities of $S$ and $C$, in other words, $\dim
\ker S^*$ and $\dim \ker C$ are unitary invariants of $T$, and the
nilpotent operator is uniquely determined by $T$ up to a
quasi-similarity. For earlier results on contractions with constant
characteristic functions see \cite{BM}, \cite{T} and \cite{WP}.

Recall that a pure isometry is a unilateral shift of some
multiplicity and a pure coisometry is the adjoint of a pure
isometry.  Recall also that a contraction $T$ on a Hilbert space
$\clh$, (i.e., $\|T h\| \leq \|h\|$ for all $h$ in $\clh$) is
completely nonunitary (c.n.u.) if there is no nontrivial reducing
subspace $\clm$ of $\clh$ for $T$ such that $T|_{\clm}$ is unitary.

In this paper we shall adopt a second approach to prove the theorem
stated in the abstract, based essentially on new factorizations of
characteristic functions of upper triangular $3 \times 3$ block
contractions (see Theorem \ref{mt-1}).

Before we continue we recall the notion of the characteristic
function of a contraction. Consider a contraction $T$ on a Hilbert
space $\clh$. The defect operators $D_T$ and $D_{T^*}$ and the
defect spaces $\cld_T$ and $\cld_{T^*}$ of $T$ are defined by
\[
D_T = (I_{\clh} - T^*T)^{\frac{1}{2}}, \quad \quad \mbox{and}
\quad\quad D_{T^*} = (I_{\clh} - T T^*)^{\frac{1}{2}},\]and
\[
\cld_T = \overline{\mbox{Ran} D_T}, \quad \quad \mbox{and} \quad
\quad \cld_{T^*} = \overline{\mbox{Ran} D_{T^*}},
\]
respectively. Then the \textit{characteristic function} of the
contraction $T$ is the $\cll(\cld_T, \cld_{T^*})$-valued contractive
analytic function defined by
\[
\Theta_T (z) = [ - T + z D_{T^*} (I_{\clh} - z T^*)^{-1}
D_{T}]|_{\cld_T} \quad \quad (z \in \mathbb{D}).
\]
In particular, $\Theta_T$ is a $\cll(\cld_T, \cld_{T^*})$-valued
bounded analytic function on $\mathbb{D}$ (see \cite{NF}). Moreover,
the characteristic function $\Theta_T$ is purely contractive, that
is,
\[
\|\Theta_T(0) \eta \| < \|\eta\|\quad \quad (\eta \in \cld_T, \eta
\neq 0).
\]

Let $\Theta \colon\  \mathbb{D} \raro \cll(\clm, \clm_*)$ and $\Psi
: \mathbb{D} \raro \cll(\cln, \cln_*)$ be two operator valued
analytic functions on $\mathbb{D}$. We say that $\Theta$ and $\Psi$
coincide and write $\Theta \cong \Psi$ if there exist two unitary
operators $\tau \colon\ \clm \raro \cln$ and $\tau_* \colon\ \clm_*
\raro \cln_*$ such that
\[\Theta(z) = \tau_*^{-1} \Psi(z) \tau \quad \quad (z \in
\mathbb{D}),\]or, equivalently, for all $z \in \mathbb{D}$ the
following diagram commutes:
\[\begin{CD}\clm @>\Theta(z) >> \clm_*\\ @ V \tau VV @V \tau_* VV\\
\cln @>\Psi (z) >>  \cln_{*}
\end{CD}\]

The characteristic function is a complete unitary invariant in the
following sense (see \cite{NF}, Theorem 3.4): Two c.n.u.
contractions $T$ on $\clh$ and $R$ on $\clk$ are unitarily
equivalent (that is, there is a unitary operator $U$ from $\clh$ to
$\clk$ such that $T = U^* R U$) if and only if
\[
\Theta_T \cong \Theta_R.
\]
Moreover, for a given $\cll(\cle, \cle_*)$-valued purely contractive
analytic function $\Theta$ defined on $\mathbb{D}$, there exists a
c.n.u. contraction $T$ on some Hilbert space, explicitly determined
by $\Theta$, such that $\Theta_T$ coincides with $\Theta$.

Contractive operator valued analytic functions play an important
role in operator theory and serve as a bridge between operator
theory and function theory in terms of systems theory and
interpolation theory (cf. \cite{DR}, \cite{NF}, \cite{NV}).

The class of nilpotent contractions yields a natural set of examples
of operators that have polynomial characteristic functions. Indeed,
let $N$ be a contraction and a nilpotent operator of order $m$, $m
\geq 1$, that is, $\|N \| \leq 1$, $N^m = 0$, and $N^{m-1} \neq 0$.
The characteristic function $\Theta_N$ of $N$ is given by
\[
\begin{split}
\Theta_N(z) & =  [- N + z D_{N^*} (I_{\clh} - z N^*)^{-1}
D_{N}]|_{\cld_N} \\
& = [- N + \sum_{p = 0}^\infty z^{p+1} D_{N^*} N^{*p}
D_N]|_{\cld_N}\\
& = [-N + \sum_{p = 0}^{m-1} z^{p+1} D_{N^*} N^{*p} D_N]|_{\cld_N},
\end{split}
\]
for all $z \in \mathbb{D}$. Therefore $\Theta_N$ is a polynomial in
$z$ of degree at most $m$ with operator coefficients.

From this viewpoint, it is important to understand, up to unitary
equivalance, the analytic structure of polynomial characteristic
functions of contractions. The main goal of the present paper is to
address this issue. More specifically, in Theorem \ref{main-thm} we
prove: If the characteristic function $\Theta_T$ of a c.n.u.
contraction $T$ is a polynomial of degree $m$, then there exist a
Hilbert space $\clm$, a nilpotent operator $N$ of order $m$, a
co-isometry $V_1 \in \cll(\cld_{N^*} \oplus \clm, \cld_{T^*})$, and
an isometry $V_2 \in \cll(\cld_{T}, \cld_N \oplus \clm)$, such that
\[
\Theta_T = V_1 \begin{bmatrix} \Theta_N & 0
\\ 0 & I_{\clm} \end{bmatrix} V_2.
\]

Along the way we prove the following factorization result for
characteristic functions (see Theorem \ref{mt-1}): Let $\clh_1,
\clh_0, \clh_{-1}$ be three Hilbert spaces and set $\clh = \clh_{1}
\oplus \clh_0 \oplus \clh_{-1}$. Let
\[
T = \begin{bmatrix} S & * & *\\0 & N & *\\0& 0& C
\end{bmatrix}
\]
be any contraction on $\clh$ with the above matricial form. Then the
characteristic function $\Theta_T$ of $T$ and
\[
\begin{bmatrix}\Theta_C & 0 \\ 0 & I_{\cle_1} \end{bmatrix} U_1
\begin{bmatrix} \Theta_N & 0
\\ 0 & I_{\clm} \end{bmatrix} U_2 \begin{bmatrix}\Theta_S & 0 \\ 0 & I_{\cle_2}
\end{bmatrix}.
\]
coincide, where $\cle_1, \cle_2$ and $\clm$ are Hilbert spaces, and
$U_1 \in \cll(\cld_{N^*} \oplus \clm, \cld_C \oplus \cle_1)$ and
$U_2 \in \cll(\cld_{S^*} \oplus \cle_2, \cld_N \oplus \clm)$ are
unitary operators.

Our results rely on the upper triangular representation of operators
with polynomial characteristic functions (see Theorem
\ref{poly-class}) and a factorization of characteristic functions of
upper triangular $2 \times 2$ block contractions due to Sz.-Nagy and
the first author (see Theorem \ref{block-cont2}).

The rest of this paper is organized as follows:. In Section 2, we
give the factorization of the characteristic function of an upper
triangular $3 \times 3$ block contraction on Hilbert space. Our main
result is given in Section 3, and provides a complete analytic
characterization of polynomial characteristic functions for c.n.u.
contractions on Hilbert space.

\newsection{Factorizations of characteristic functions}

We start by recalling some known facts about upper triangular $2
\times 2$ block contractions, since they will be frequently used in
what follows.

The first is a classification of $2 \times 2$ block contractions.
This is the content of Theorem 1 in \cite{NF1} (also see Chapter IV,
Lemma 2.1 in \cite{FF}).

\begin{Theorem}\label{block-cont}
Let $\clh_1$ and $\clh_2$ be Hilbert spaces and let $T =
\begin{bmatrix}T_1 & X\\ 0 & T_2\end{bmatrix}$ be a bounded linear operator on $\clh_1 \oplus \clh_2$.
Then $T$ is a contraction if and only if $T_1$ and $T_2$ are
contractions and \[X = D_{T_1^*} \Gamma D_{T_2},\]for some
contraction $\Gamma$ from $\cld_{T_2}$ to $\cld_{T_1^*}$.
\end{Theorem}

The second key tool used in our development is the factorization of
characteristic functions of $2 \times 2$ block contractions (see
Theorem 2 in \cite{NF1}).

\begin{Theorem}\label{block-cont2}
Let $\clh_1$ and $\clh_2$ be Hilbert spaces, let $T =
\begin{bmatrix}T_1 & X\\ 0 & T_2\end{bmatrix}$ be a contraction on $\clh_1 \oplus
\clh_2$, and let $X = D_{T_1^*} \Gamma D_{T_2}$ for some contraction
$\Gamma \in \cll(\cld_{T_1^*}, \cld_{T_2})$. Then there exist
unitary operators $\tau \in \cll(\cld_T, \cld_{T_1} \oplus
\cld_{\Gamma})$ and $\tau_* \in \cll(\cld_{T^*}, \cld_{T_2^*} \oplus
\cld_{\Gamma^*})$ such that
\[\Theta_T(z) = \tau_*^{-1} \begin{bmatrix} \Theta_{T_2}(z) & 0 \\ 0
& I_{\cld_{\Gamma^*}}\end{bmatrix} J[\Gamma] \begin{bmatrix}
\Theta_{T_1}(z) & 0 \\ 0 & I_{\cld_{\Gamma}}\end{bmatrix} \tau \quad
\quad (z \in \mathbb{D}),\]where \[J[\Gamma] = \begin{bmatrix}
\Gamma^* & D_{\Gamma}
\\ D_{\Gamma^*} & - \Gamma\end{bmatrix} \in \cll(\cld_{T_1^*} \oplus \cld_{\Gamma}, \cld_{T_2} \oplus \cld_{\Gamma^*}).\]
\end{Theorem}

Recall that if $A$ is a contraction from $\clh$ to $\clk$ then
\begin{equation}\label{JH}
J[A] = \begin{bmatrix} A^* & D_{A}\\ D_{A^*} & - A
\end{bmatrix}
\end{equation}
is a unitary operator from $\clk \oplus \cld_{A}$ to $\clh \oplus
\cld_{A^*}$ (see Halmos \cite{H}).

We are now ready to prove our first factorization result.

\begin{Theorem}\label{mt-1}
Let $\clh_1, \clh_0, \clh_{-1}$ be Hilbert spaces, and let $\clh =
\clh_{1} \oplus \clh_0 \oplus \clh_{-1}$. Let
\[
T = \begin{bmatrix} S & * & *\\0 & N & *\\0& 0& C
\end{bmatrix},
\]
be a contraction on $\clh$. Then there exist three Hilbert spaces
$\cle_1, \cle_2$ and $\clm$ and two unitary operators $U_1 \in
\cll(\cld_{N^*} \oplus \clm, \cld_C \oplus \cle_1)$ and $U_2 \in
\cll(\cld_{S^*} \oplus \cle_2, \cld_N \oplus \clm)$ such that
\[\Theta_T \cong  \begin{bmatrix}\Theta_C & 0 \\ 0 & I_{\cle_1} \end{bmatrix} U_1
\begin{bmatrix} \Theta_N & 0
\\ 0 & I_{\clm} \end{bmatrix} U_2 \begin{bmatrix}\Theta_S & 0 \\ 0 & I_{\cle_2}
\end{bmatrix}.\]
\end{Theorem}
\NI\textsf{Proof.} Set

\[
\clk_1 = \clh_1 \oplus \clh_0,
\]
\[
T =\begin{bmatrix} T_1 & X_1
\\ 0 & C\end{bmatrix} \in \cll(\clk_1 \oplus \clh_{-1}),\]
and
\[T_1 = \begin{bmatrix} S & X\\ 0 & N\end{bmatrix} \in \cll(\clh_1 \oplus \clh_0),\]
where $X_1 \in \cll(\clh_{-1}, \clk_1)$ and $X \in \cll(\clh_0,
\clh_1)$. Theorem \ref{block-cont} implies that there exist
contractions $\Gamma_1 \in \cll(\cld_C, \cld_{T_1^*})$ and $\Gamma
\in \cll(\cld_N, \cld_{S^*})$ such that
\[X_1 = D_{T_1^*} \Gamma_1 D_C,\]and \[X = D_{S^*} \Gamma D_N.\]
By Theorem \ref{block-cont2} there exist unitary operators
\begin{equation}\label{tau1}\tau_1 : \cld_T \raro \cld_{T_1} \oplus \cld_{\Gamma_1}, \quad
\tau_{1*} : \cld_{T^*} \raro \cld_{C^*} \oplus
\cld_{\Gamma_1^*},\end{equation} and
\begin{equation}\label{tau}\tau : \cld_{T_1} \raro \cld_{S} \oplus \cld_{\Gamma}, \quad
\tau_{*} : \cld_{T_1^*} \raro \cld_{N^*} \oplus
\cld_{\Gamma^*},\end{equation} such that
\[\Theta_T(z) = \tau^{-1}_{1*} \begin{bmatrix} \Theta_C(z) & 0 \\ 0 &
I_{\cld_{\Gamma_1^*}}\end{bmatrix} J[\Gamma_1] \begin{bmatrix}
\Theta_{T_1}(z) & 0
\\ 0 & I_{\cld_{\Gamma_1}}\end{bmatrix} \tau_1,\]
and
\[\Theta_{T_1}(z) = \tau^{-1}_{*} \begin{bmatrix} \Theta_N(z) & 0 \\ 0 &
I_{\cld_{\Gamma^*}}\end{bmatrix} J[\Gamma] \begin{bmatrix}
\Theta_S(z) & 0
\\ 0 & I_{\cld_{\Gamma}}\end{bmatrix} \tau,\]
for all $z \in \mathbb{D}$, where
\[J[\Gamma_1] = \begin{bmatrix}
\Gamma_1^* & D_{\Gamma_1}
\\ D_{\Gamma_1^*} & - \Gamma_1\end{bmatrix} \in \cll(\cld_{T_1^*} \oplus \cld_{\Gamma_1}, \cld_{C} \oplus \cld_{\Gamma_1^*}),\]
and
\[J[\Gamma] = \begin{bmatrix}
\Gamma^* & D_{\Gamma}
\\ D_{\Gamma^*} & - \Gamma\end{bmatrix} \in \cll(\cld_{S^*} \oplus \cld_{\Gamma}, \cld_N \oplus \cld_{\Gamma^*}),\]
are unitary operators (see (\ref{JH})). Now setting
\[
\Phi_S(z) =
\begin{bmatrix}
\Theta_S(z) & 0
\\ 0 & I_{\cld_{\Gamma}}\end{bmatrix},
\Phi_N(z) = \begin{bmatrix} \Theta_N(z) & 0 \\ 0 &
I_{\cld_{\Gamma^*}}\end{bmatrix}, \mbox{~ and ~ } \Phi_C(z) =
\begin{bmatrix} \Theta_C(z) & 0 \\ 0 &
I_{\cld_{\Gamma_1^*}}\end{bmatrix},
\]
for all $z \in \mathbb{D}$, we get
\[
\begin{split}
\Theta_T(z) & =  \tau^{-1}_{1*} \Phi_C(z) J[\Gamma_1]
\begin{bmatrix} \Theta_{T_1}(z) & 0
\\ 0 & I_{\cld_{\Gamma_1}}\end{bmatrix} \tau_1 \\ &
= \tau^{-1}_{1*} \Phi_C(z) J[\Gamma_1] \begin{bmatrix} \tau^{-1}_{*}
\Phi_N(z) J[\Gamma] \Phi_S(z) \tau & 0
\\ 0 & I_{\cld_{\Gamma_1}}\end{bmatrix} \tau_1 \\ &
= \tau^{-1}_{1*} \Phi_C(z) \Big( J[\Gamma_1]
\begin{bmatrix}\tau^{-1}_* & 0 \\ 0 & I_{\cld_{\Gamma_1}} \end{bmatrix} \Big) \begin{bmatrix}
\Phi_N(z) J[\Gamma] \Phi_S(z) \tau & 0
\\ 0 & I_{\cld_{\Gamma_1}}\end{bmatrix} \tau_1 \\ &
= \tau^{-1}_{1*} \Phi_C(z)  U_1 \begin{bmatrix} \Phi_N(z) J[\Gamma]
\Phi_S(z)\tau & 0
\\ 0 & I_{\cld_{\Gamma_1}}\end{bmatrix} \tau_1,
\end{split}
\]
for all $z \in \mathbb{D}$, where $U_1 \in \cll((\cld_{N^*} \oplus
\cld_{\Gamma^*}) \oplus \cld_{\Gamma_1}, \cld_{C} \oplus
\cld_{\Gamma_1^*})$ is the unitary operator defined by

\[
U_1 = J[\Gamma_1]
\begin{bmatrix}\tau^{-1}_* & 0 \\ 0 & I_{\cld_{\Gamma_1}}
\end{bmatrix}.
\]
Hence
\[
\begin{split}
\Theta_T(z) & = \tau^{-1}_{1*} \Phi_C(z) U_1 \begin{bmatrix}
\Phi_N(z) J[\Gamma] \Phi_S(z) \tau & 0
\\ 0 & I_{\cld_{\Gamma_1}}\end{bmatrix} \tau_1
\\ &
= \tau^{-1}_{1*} \Phi_C(z) U_1 \begin{bmatrix} \Phi_N(z) & 0 \\ 0 &
I_{\cld_{\Gamma_1}} \end{bmatrix} \begin{bmatrix}J[\Gamma] & 0 \\
0 & I_{\cld_{\Gamma_1}}
\end{bmatrix}
\begin{bmatrix} \Phi_S(z) & 0 \\ 0 & I_{\cld_{\Gamma_1}} \end{bmatrix}
\begin{bmatrix} \tau & 0 \\ 0 & I_{\cld_{\Gamma_1}}\end{bmatrix} \tau_1
\\ &
= \tau^{-1}_{1*} \Phi_C(z) U_1 \begin{bmatrix} \Theta_N(z) & 0 \\ 0
&
I_{\cld_{\Gamma^*} \oplus \cld_{\Gamma_1}} \end{bmatrix} \begin{bmatrix}J[\Gamma] & 0 \\
0 & I_{\cld_{\Gamma_1}}
\end{bmatrix}
\begin{bmatrix} \Theta_S(z) & 0 \\ 0 & I_{\cld_{\Gamma} \oplus
\cld_{\Gamma_1}} \end{bmatrix}
\begin{bmatrix} \tau & 0 \\ 0 & I_{\cld_{\Gamma_1}}\end{bmatrix} \tau_1,
\end{split}
\]
for all $z \in \mathbb{D}$. Let $U_2 \in \cll((\cld_{S^*} \oplus
\cld_{\Gamma}) \oplus \cld_{\Gamma_1}, (\cld_{N} \oplus
\cld_{\Gamma^*}) \oplus \cld_{\Gamma_1})$ and $\tilde{\tau}_1 \in
\cll(\cld_T, (\cld_S \oplus \cld_{\Gamma}) \oplus \cld_{\Gamma_1})$
be unitary operators defined by
\[
U_2 = \begin{bmatrix}J[\Gamma] & 0 \\
0 & I_{\cld_{\Gamma_1}}
\end{bmatrix},
\]
and
\[
\tilde{\tau}_1 =
\begin{bmatrix}\tau & 0 \\ 0 & I_{\cld_{\Gamma_1}}
\end{bmatrix} \tau_1,
\]
respectively. Hence we obtain
\begin{equation}\label{eq-TT}
\Theta_T(z) = \tau^{-1}_{1*} \Big( \begin{bmatrix} \Theta_C(z) & 0
\\ 0 & I_{\cld_{\Gamma_1^*}} \end{bmatrix} U_1 \begin{bmatrix} \Theta_N(z) & 0 \\ 0 &
I_{\cld_{\Gamma^*} \oplus \cld_{\Gamma_1}} \end{bmatrix} U_2
\begin{bmatrix} \Theta_S(z) & 0 \\ 0 & I_{\cld_{\Gamma} \oplus
\cld_{\Gamma_1}} \end{bmatrix} \Big) \tilde{\tau}_1,
\end{equation}
for all $z \in \mathbb{D}$, and therefore
\[
\Theta_T \cong \begin{bmatrix} \Theta_C & 0
\\ 0 & I_{\cld_{\Gamma_1^*}} \end{bmatrix} U_1 \begin{bmatrix} \Theta_N & 0 \\ 0 &
I_{\cld_{\Gamma^*} \oplus \cld_{\Gamma_1}} \end{bmatrix} U_2
\begin{bmatrix} \Theta_S & 0 \\ 0 & I_{\cld_{\Gamma} \oplus
\cld_{\Gamma_1}} \end{bmatrix},
\]
holds. Setting $\cle_1 = \cld_{\Gamma_1^*}$, $\clm = \cld_{\Gamma^*}
\oplus \cld_{\Gamma_1}$ and $\cle_2 = \cld_{\Gamma} \oplus
\cld_{\Gamma_1}$ in the above, we conclude the proof of the theorem.
\qed

Of particular interest is the case when $S$ and $C^*$ are pure
isometries.

\begin{Corollary}\label{Cor-mt}
With the hypotheses of Theorem \ref{mt-1}, let also assume that $S$
and $C^*$ are pure isometries. Then there exist a Hilbert space
$\clm$, a co-isometry $V_1 \in \cll(\cld_{N^*} \oplus \clm,
\cld_{T^*})$, and an isometry $V_2 \in \cll(\cld_{T}, \cld_N \oplus
\clm)$, such that
\[
\Theta_T = V_1
\begin{bmatrix} \Theta_N & 0
\\ 0 & I_{\clm} \end{bmatrix} V_2.
\]
\end{Corollary}
\NI\textsf{Proof.} Notice that since $\cld_{C^*} =
\{0_{\clh_{-1}}\}$ and $\cld_{S} = \{0_{\clh_1}\}$, the
characteristic functions $\Theta_C : \mathbb{D} \raro \cll(\cld_C,
\cld_{C^*})$ of $C$ and $\Theta_S : \mathbb{D} \raro \cll(\cld_S,
\cld_{S^*})$ of $S$ are identically zero, that is,
\[
0_C := \Theta_C \equiv 0 : \cld_C \raro \{0_{\clh_{-1}}\},
\mbox{~and~} 0_S := \Theta_S \equiv 0 : \{0_{\clh_1}\} \raro
\cld_{S^*}.
\]
Furthermore, the unitary operators in (\ref{tau1}) and (\ref{tau})
become
\begin{equation}\label{tau1-II}
\tau_1 : \cld_T \raro \cld_{T_1} \oplus \cld_{\Gamma_1}, \quad
\tau_{1*} : \cld_{T^*} \raro \{0_{\clh_{-1}}\} \oplus
\cld_{\Gamma_1^*},\end{equation} and
\begin{equation}\label{tau-II}\tau : \cld_{T_1} \raro \{0_{\clh_1}\} \oplus \cld_{\Gamma}, \quad
\tau_{*} : \cld_{T_1^*} \raro \cld_{N^*} \oplus
\cld_{\Gamma^*}.
\end{equation}
This along with (\ref{eq-TT}) yields
\[
\begin{split}\Theta_T & =
\tau_{1*}^{-1} \begin{bmatrix} 0_C & 0
\\ 0 & I_{\cld_{\Gamma_1^*}} \end{bmatrix} U_1
\begin{bmatrix} \Theta_N & 0 \\ 0 &
I_{\cld_{\Gamma^*} \oplus \cld_{\Gamma_1}} \end{bmatrix} U_2
\begin{bmatrix} 0_S & 0 \\ 0 & I_{\cld_{\Gamma} \oplus
\cld_{\Gamma_1}} \end{bmatrix} \tilde{\tau}_1
\\& =
V_1 \begin{bmatrix} \Theta_N & 0 \\ 0 & I_{\clm} \end{bmatrix} V_2,
\end{split}\]
where
\[
\clm = \cld_{\Gamma^*} \oplus \cld_{\Gamma_1},
\]
\[
V_1 = \tau_{1*}^{-1} \begin{bmatrix} 0_C & 0
\\ 0 & I_{\cld_{\Gamma_1^*}} \end{bmatrix} U_1
\in \cll((\cld_{N^*} \oplus \cld_{\Gamma^*}) \oplus \cld_{\Gamma_1},
\cld_{T^*})
\]
and
\[
V_2 =  U_2
\begin{bmatrix} 0_S & 0 \\ 0 & I_{\cld_{\Gamma} \oplus
\cld_{\Gamma_1}} \end{bmatrix} \tilde{\tau}_1 \in \cll(\cld_{T},
(\cld_N \oplus \cld_{\Gamma^*}) \oplus \cld_{\Gamma_1}).
\]
Now using $0_C 0_C^* = I_{\cld_{C^*}} = I_{\{0_{\clh_{-1}}\}}$ and
$0_S^* 0_S = I_{\cld_S} = I_{\{0_{\clh_{1}}\}}$ along with
(\ref{tau1-II}) and (\ref{tau-II}) we readily see that $V_1 V_1^* =
I_{\cld_{T^*}}$ and $V_2^* V_2 = I_{\cld_T}$. This completes the
proof of the corollary. \qed

\newsection{Polynomial characteristic functions}

For the readers convenience, we first state the main result of
\cite{FS}.

\begin{Theorem}\label{poly-class}
Let $T$ be a c.n.u. contraction on a Hilbert space $\clh$. Then the
characteristic function $\Theta_T$ of $T$ is a polynomial of degree
$m$ if and only if there exist three closed subspaces $\clh_1,
\clh_0, \clh_{-1}$ of $\clh$ with $\clh = \clh_1 \oplus \clh_0
\oplus \clh_{-1}$, a pure isometry $S$ in $\cll(\clh_1)$, a
nilpotent $N$ of order $m$ in $\cll(\clh_0)$, and a pure co-isometry
$C$ in $\cll(\clh_{-1})$, such that $T$ has the matrix
representation
\begin{equation*}
T = \begin{bmatrix} S & * & *\\0 & N & *\\0& 0& C
\end{bmatrix}.
\end{equation*}
\end{Theorem}

We are now ready for the main theorem on analytic description of
contractions which have polynomial characteristic functions.

\begin{Theorem}\label{main-thm}
Let $T$ be a c.n.u. contraction on a Hilbert space $\clh$. If the
characteristic function $\Theta_T$ of $T$ is a polynomial of degree
$m$, then there exist a Hilbert space $\clm$, a nilpotent operator
$N$ of order $m$, a co-isometry $V_1 \in \cll(\cld_{N^*} \oplus
\clm, \cld_{T^*})$, and an isometry $V_2 \in \cll(\cld_{T}, \cld_N
\oplus \clm)$, such that
\[
\Theta_T = V_1 \begin{bmatrix} \Theta_N & 0
\\ 0 & I_{\clm} \end{bmatrix} V_2.
\]
\end{Theorem}
\NI\textsf{Proof.} Let $T$ be a c.n.u. contraction such that the
characteristic function $\Theta_T$ of $T$ is a polynomial of degree
$m$. According to Theorem \ref{poly-class} there exist closed
subspaces $\clh_1, \clh_0, \clh_{-1}$ of $\clh$ such that $\clh =
\clh_1 \oplus \clh_0 \oplus \clh_{-1}$ and such that with respect to
that decomposition, $T$ admits the matrix representation
\[
T = \begin{bmatrix} S & * & *\\0 & N & *\\0& 0& C
\end{bmatrix},
\]
where $S \in \cll(\clh_1)$ is a  pure isometry, $N \in \cll(\clh_0)$
is a nilpotent operator of order $m$, and  $C \in \cll(\clh_{-1})$
is a pure coisometry. The result now follows from Corollary
\ref{Cor-mt}. \qed

\begin{Remark}\label{remark-3.2}
The converse of the above theorem is not true in full generality:
Let $\clm$ be an infinite dimensional separable Hilbert space, and
let $T$ be a c.n.u. contraction with infinite dimensional defect
spaces (for example, one can consider $T = S \oplus S^*$ on
$H^2_{\cle}(\mathbb{D}) \oplus H^2_{\cle}(\mathbb{D})$, where $\cle$
is an infinite dimensional Hilbert space, $H^2_{\cle}(\mathbb{D})$
is the $\cle$-valued Hardy space, and $S$ is the shift operator on
$H^2_{\cle}(\mathbb{D})$). Let $N$ be a nilpotent operator of order
$m$ and let
\[
V_2 = \begin{bmatrix} V_{21}\\V_{22}\end{bmatrix} : \cld_T \raro
\cld_N \oplus \clm,
\]
be an isometry, where
\[
\|V_{21} \eta \| = \| V_{22} \eta \| \quad \quad (\eta \in \cld_T).
\]
Also, let $V_1 : \cld_{N^*} \oplus \clm \raro \cld_{T^*}$ be a
coisometry with $\ker V_1 = \cld_{N^*}$. If
\[
\Theta_T = V_1 \begin{bmatrix} \Theta_N & 0
\\ 0 & I_{\clm} \end{bmatrix} V_2,
\]
then $\Theta_T$ is a polynomial of degree $0$.

\noindent However, it is easy to see that the following weak
converse of Theorem \ref{main-thm} is true: Let $T$ be a c.n.u.
contraction on a Hilbert space $\clh$. Let
\[
\Theta_T = V_1 \begin{bmatrix} \Theta_N & 0
\\ 0 & I_{\clm} \end{bmatrix} V_2,
\]
for some Hilbert space $\clm$, nilpotent operator $N$ of order $m$,
a co-isometry $V_1 \in \cll(\cld_{N^*} \oplus \clm, \cld_{T^*})$,
and an isometry $V_2 \in \cll(\cld_{T}, \cld_N \oplus \clm)$. Then
the characteristic function $\Theta_T$ of $T$ is a polynomial of
degree less than or equal to $m$.
\end{Remark}

It is important to note that the conclusion of Theorem
\ref{main-thm} depends explicitly on the decomposition of $T$ as
used in the proof of Theorem \ref{mt-1}. With the same setting as in
Theorem \ref{main-thm}, below we will show that the same conclusion
holds for the following decomposition of $T$:
\[
T = \begin{bmatrix} S & X_{-1}
\\ 0 & T_{-1} \end{bmatrix} =
\begin{bmatrix} S & D_{S^*} \Gamma_{-1} D_{T_{-1}}
\\ 0 & T_{-1} \end{bmatrix} \in \cll(\clh_0 \oplus \clk_{-1}),\]
where $\clk_{-1} = \clh_0 \oplus \clh_{-1}$,
\[
T_{-1} = \begin{bmatrix} N & X\\ 0 & C\end{bmatrix} =
\begin{bmatrix} N & D_{N^*} \Gamma D_C \\ 0 & C\end{bmatrix} \in
\cll(\clh_1 \oplus \clh_0),
\]
and $X_{-1} = D_{S^*} \Gamma_{-1} D_{T_1}$, $X = D_{N^*} \Gamma
D_C$, and $\Gamma_{-1}$ in $\cll(\cld_{T_{-1}}, \cld_{S^*})$ and
$\Gamma$ in $\cll(\cld_C, \cld_{N^*})$ are a pair of contractions.
In this case, again by Theorem \ref{block-cont2}, we have
\begin{equation}\label{T}
\Theta_T = \tau^{-1}_{-1*}
\begin{bmatrix} \Theta_{T_{-1}} & 0 \\ 0 &
I_{\cld_{\Gamma_{-1}^*}}\end{bmatrix} J[\Gamma_{-1}] \begin{bmatrix}
0_S & 0
\\ 0 & I_{\cld_{\Gamma_{-1}}}\end{bmatrix} \tau_{-1},
\end{equation}
and
\begin{equation}\label{T-1}
\Theta_{T_{-1}} = \tau^{-1}_{*} \begin{bmatrix} 0_C & 0 \\
0 & I_{\cld_{\Gamma^*}}\end{bmatrix} J[\Gamma] \begin{bmatrix}
\Theta_N & 0
\\ 0 & I_{\cld_{\Gamma}}\end{bmatrix} \tau,
\end{equation}
where
\begin{equation}\label{tau1-III}\tau_{-1} : \cld_T \raro \{0_{\clh_1}\} \oplus \cld_{\Gamma_{-1}}, \quad
\tau_{-1*} : \cld_{T^*} \raro \cld_{T_{-1}^*} \oplus
\cld_{\Gamma_{-1}^*},\end{equation} and
\begin{equation}\label{tau-III}\tau : \cld_{T_{-1}} \raro \cld_N \oplus \cld_{\Gamma}, \quad
\tau_{*} : \cld_{T_{-1}^*} \raro \{0_{\clh_{-1}}\} \oplus
\cld_{\Gamma^*},\end{equation}
are unitary operators. Moreover
\[
J[\Gamma_{-1}] = \begin{bmatrix}
\Gamma_{-1}^* & D_{\Gamma_{-1}}
\\ D_{\Gamma_{-1}^*} & - \Gamma_{-1}\end{bmatrix} \in \cll(\cld_{S^*}
\oplus \cld_{\Gamma_{-1}}, \cld_{T_{-1}} \oplus \cld_{\Gamma_{-1}^*}),
\]
and
\[
J[\Gamma] = \begin{bmatrix}
\Gamma^* & D_{\Gamma}
\\ D_{\Gamma^*} & - \Gamma\end{bmatrix} \in \cll(\cld_{N^*} \oplus \cld_{\Gamma}, \cld_C \oplus \cld_{\Gamma^*}).
\]
By setting
\[
\Psi_0 = \begin{bmatrix} 0_C & 0
\\ 0 & I_{\cld_{\Gamma^*}}\end{bmatrix}
\mbox{~and~} \Psi_N = \begin{bmatrix} \Theta_N & 0
\\ 0 & I_{\cld_{\Gamma}}\end{bmatrix},
\]
and using (\ref{T}) and (\ref{T-1}) we obtain
\[
\begin{split}
\Theta_T & = \tau^{-1}_{-1*} \begin{bmatrix} \Theta_{T_{-1}} & 0 \\
0 & I_{\cld_{\Gamma_{-1}^*}}\end{bmatrix} J[\Gamma_{-1}]
\begin{bmatrix} 0_S & 0
\\ 0 & I_{\cld_{\Gamma_{-1}}}\end{bmatrix} \tau_{-1} \\
& =
\tau^{-1}_{-1*} \begin{bmatrix} \tau^{-1}_{*} \Psi_0 J[\Gamma]
\Psi_N \tau & 0
\\ 0 & I_{\cld_{\Gamma_{-1}^*}}\end{bmatrix} J[\Gamma_{-1}] \begin{bmatrix}
0_S & 0
\\ 0 & I_{\cld_{\Gamma_{-1}}}\end{bmatrix} \tau_{-1}
\\ &
= \tau^{-1}_{-1*} \begin{bmatrix} \tau_*^{-1} & 0 \\ 0 &
I_{\cld_{\Gamma_{-1}^*}}\end{bmatrix}
\begin{bmatrix} \Psi_0 & 0 \\ 0 & I_{\cld_{\Gamma_{-1}^*}}
\end{bmatrix}
\begin{bmatrix} J[\Gamma] & 0 \\ 0 &
I_{\cld_{\Gamma_{-1}^*}}\end{bmatrix}
\begin{bmatrix} \Psi_N & 0 \\ 0 & I_{\cld_{\Gamma_{-1}^*}}
\end{bmatrix}
\begin{bmatrix} \tau & 0 \\ 0 & I_{\cld_{\Gamma_{-1}^*}}
\end{bmatrix}
 J[\Gamma_{-1}] \begin{bmatrix}
0_S & 0
\\ 0 & I_{\cld_{\Gamma_{-1}}}\end{bmatrix} \tau_{-1}
\\ &
= \tilde{V}_1 \begin{bmatrix} \Psi_N & 0 \\ 0 &
I_{\cld_{\Gamma_{-1}^*}}
\end{bmatrix} \tilde{V}_2
\\ &
= \tilde{V}_1 \begin{bmatrix} \Theta_N & 0 \\ 0 & I_{\cld_{\Gamma}
\oplus \cld_{\Gamma_{-1}^*}}
\end{bmatrix} \tilde{V}_2,
\end{split}
\]
where
\[
\begin{split}
\tilde{V}_1 & = \tau^{-1}_{-1*} \begin{bmatrix} \tau_*^{-1} & 0 \\ 0
& I_{\cld_{\Gamma_{-1}^*}}\end{bmatrix}
\begin{bmatrix} \Psi_0 & 0 \\ 0 & I_{\cld_{\Gamma_{-1}^*}}
\end{bmatrix}
\begin{bmatrix} J[\Gamma] & 0 \\ 0 &
I_{\cld_{\Gamma_{-1}^*}}\end{bmatrix}
\\
& = \tau^{-1}_{-1*} \begin{bmatrix} \tau_*^{-1} & 0 \\
0 & I_{\cld_{\Gamma_{-1}^*}}\end{bmatrix}
\begin{bmatrix}0_C & 0 \\ 0 & I_{\cld_{\Gamma^*} \oplus \cld_{\Gamma_{-1}^*}}
\end{bmatrix}
\begin{bmatrix} J[\Gamma] & 0 \\ 0 &
I_{\cld_{\Gamma_{-1}^*}}\end{bmatrix},
\end{split}
\]
and
\[
\tilde{V}_2 = \begin{bmatrix} \tau & 0 \\ 0 &
I_{\cld_{\Gamma_{-1}^*}}\end{bmatrix} J[\Gamma_{-1}]
\begin{bmatrix} 0_S & 0 \\ 0 & I_{\cld_{\Gamma_{-1}}} \end{bmatrix} \tau_{-1}.
\]
Hence
\[
\begin{split}
\Theta_T & = \tilde{V}_1 \begin{bmatrix} \Theta_N & 0 \\
0 & I_{\cld_{\Gamma} \oplus \cld_{\Gamma_{-1}^*}} \end{bmatrix}
\tilde{V}_2
\\
& =
\tilde{V}_1 \begin{bmatrix} \Theta_N & 0 \\
0 & I_{\tilde{\clm}}
\end{bmatrix} \tilde{V}_2
\end{split},
\]
where $\tilde{\clm} = \cld_{\Gamma} \oplus \cld_{\Gamma_{-1}^*}$.
Finally, by virtue of (\ref{tau1-III}) and (\ref{tau-III}), we have
that $\tilde{V}_1^*$ and $\tilde{V}_2$ are isometric operators, that
is, $\tilde{V}_1 \tilde{V}_1^* = I_{\cld_{T^*}}$ and $\tilde{V}_2^*
\tilde{V}_2 = I_{\cld_T}$. Yet, we do not know if
\[
\mbox{dim~} \tilde{\clm} = \mbox{dim~} \clm,
\]
where $\clm$ is as in Theorem \ref{main-thm}.

\vspace{0.3in}

\NI\textit{Acknowledgement:} The authors are grateful to the referee
for pointing out that the formulation of Theorem \ref{main-thm} in
the submitted version was inappropriate. The example in Remark
\ref{remark-3.2} is due to the referee.

\noindent The third author is supported in part by NBHM (National
Board of Higher Mathematics, India) grant NBHM/R.P.64/2014. The
third author is also grateful for hospitality of Texas A\&M
University, USA, during July 2016 and July-August 2015.

\end{document}